\documentclass{egpubl}
\usepackage{wait2015}

\ConferencePaper      %

\title%
{On the smallest parallel quadrangulation with minimum degree 3}

\usepackage[utf8]{inputenc}  %

\author[Kápolnai et al]
       {Richárd Kápolnai,$^1$ %
        Gábor Domokos$^2$ and
        Imre Szeberényi$^1$
        \\
        $^1$ Dept.\ of Control Engineering and Information Technology, Budapest University of Technology and Economics\\
        $^2$ Dept.\ of Mechanics, Materials \& Structures, Budapest University of Technology and Economics\\
Corresponding email: \url{kapolnai@iit.bme.hu}
       }

\usepackage{graphicx}

\PrintedOrElectronic

\usepackage{t1enc,dfadobe}

        \let\tt=\ttfamily
          
\let\bf=\bfseries

\usepackage[hidelinks]{hyperref}
\usepackage[inline]{enumitem}
\usepackage{mathrsfs} %
\usepackage{subcaption}

\newcounter{tetelc}
\newcommand{\tetel}[1]{%
\refstepcounter{tetelc}%
{\setlength{\parindent}{0cm}%
{\bf Theorem~\arabic{tetelc}.} \textit{#1}%
}}

\usepackage{acronym}
\acrodef{MUQ}[MQ]{multiquadrangulation} %

\begin{document}

\maketitle

\begin{abstract}
The identity of the smallest quadrangulation with minimum degree 3 also containing parallel edges is unknown.
However, it has already been determined that its order (the number of vertices) is between 11 and 14.
This paper narrows this domain by showing that the order is at least 12.

\end{abstract}

\section{Introduction}

\begin{figure}%
\newcommand{\qscaleratio}{0.85}
\centering
\subcaptionbox{$P_2$\label{fig:smallests:p2}}{
  \includegraphics[scale=\qscaleratio]{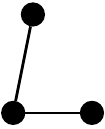}
}
\quad
\subcaptionbox{$C_4$\label{fig:smallests:square}}{
  \includegraphics[scale=\qscaleratio]{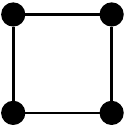}
}
\quad
\subcaptionbox{$Q_3$\label{fig:smallests:q3}}{
  \includegraphics[scale=\qscaleratio]{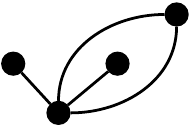}
}
\quad
\subcaptionbox{$Q_4$\label{fig:smallests:q4}}{
  \includegraphics[scale=\qscaleratio]{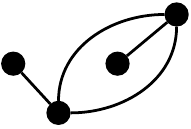}
}
\caption[The smallest quadrangulations]{\acp{MUQ} of smallest order}
\label{fig:smallests}
\end{figure}

A \emph{plane graph} is a graph whose vertices are drawn points and edges are arcs on the two dimensional plane such that no two edges meet in a point other than a common endpoint\cite{Diestel2005}.
The edges divide the planar surface into regions called \emph{faces}.
A \emph{walk} of length $l$ is a sequence of $l$ adjacent edges, and the walk is \emph{closed} if it ends in the starting vertex.
A \emph{plane quadrangulation} (or shortly \emph{quadrangulation}) is a loopless, connected, finite plane graph having every face bounded by a closed walk of length 4.
A quadrangulation without parallel edges and without repeated edges on the quadrilateral boundary walks is called a \emph{simple quadrangulation}.
We allow parallel edges, and the boundary walk may repeat edges or vertices.
If we want to emphasize that a quadrangulation may not be simple, it is called a \emph{multiquadrangulation}, abbreviated as \acs{MUQ}.
The \acp{MUQ} of smallest order are shown in \autoref{fig:smallests}.

This paper investigates the smallest \acp{MUQ} with minimum degree 3 which has parallel edges.
By \emph{smallest}, we mean smallest \emph{order}, i.e.\ the minimum number of vertices.
\autoref{tab:illus} illustrates these concepts on the graphs of \autoref{fig:smallests}.
Note that though $P_2$ is not parallel, it is not simple either according to the definition above.
This paper shows the following:

\tetel{
\label{theo:nosmallparallel}
Every \ac{MUQ} with minimum degree 3 containing parallel edges has at least 12 vertices.
}

\begin{table}
  \centering
  \begin{tabular}{p{0.15\linewidth}p{0.15\linewidth}p{0.15\linewidth}p{0.2\linewidth}}
    graph & order & parallel & min.\ degree\\
\hline
    $P_2$ & 3 & no & 1\\
    $C_4$ & 4 & no & 2\\
    $Q_3$ & 4 & yes & 1\\
    $Q_4$ & 4 & yes & 1\\
  \end{tabular}
\caption{Illustration of concepts}
\label{tab:illus}
\end{table}

The importance of \acp{MUQ} with minimum degree 3
is related to a 
mechanical classification system for convex, homogeneous bodies
introduced recently\cite{Domokos2006,Domokos2012arxiv}.
In this system, each body is mapped into its \emph{secondary equilibrium class} determined by the topology of the equilibrium points of the surface.
Such a topology is defined by a vertex-coloured \ac{MUQ}.
There are particular secondary equilibrium classes called \emph{irreducible ancestors}, from which every other secondary equilibrium class can be generated with specific transformations (detailed in other works\cite{Domokos2012arxiv,kapolnai2012periodica}).
Such an ancestor class either corresponds to a \ac{MUQ} with minimum degree 3, %
or to the \ac{MUQ} denoted by $P_2$ on \autoref{fig:smallests}, latter representing the class of the mono-monostatic body called Gömböc\cite{kapolnai2012periodica}.
For this reason, we say a \ac{MUQ} is an \emph{irreducible \ac{MUQ}} (or shortly \emph{irreducible}), if it is isomorphic to $P_2$ or its minimum degree is 3.
If an irreducible contains parallel edges, it is called a \emph{parallel irreducible}.

There are efficient methods to exhaustively enumerate simple irreducibles\cite{Nakamoto1999,Brinkmann2005}.
A highly tuned implementation called \emph{Plantri} is also available\cite{Brinkmann2007}.
However, these enumerations are incomplete as they ignore the ones containing parallel edges.
Identifying the smallest parallel irreducible would specify the limit of these incomplete methods regarding the enumeration of all \acp{MUQ}.

We mention that some pieces of the text of this paper can also be found in the submitted dissertation of the first author.

\section{Related Work}

It has already been determined that the order of the smallest parallel irreducible is between 11 and 14, detailed in this section.
First, every irreducible has at least 8 vertices\cite{Batagelj1989,kapolnai2012periodica}:

\tetel{%
  \label{prop:Batageljgen}
  The order of an irreducible \ac{MUQ} is at least 8.
}

This statement was originally made by Batagelj\cite{Batagelj1989} for simple irreducibles, and was later generalized by others\cite{kapolnai2012periodica} to include parallel irreducibles as well.
There also exists a primitive implementation\cite{kapolnai2012periodica} to enumerate every \ac{MUQ} of order at most 10.
Observing the generated data\cite{kapolnai2012periodica}, we arrive at

\tetel{%
Having the order between 8 and 10,
there are two non-isomorphic irreducible \acp{MUQ}.%
}

\begin{figure}
\centering
\subcaptionbox{}{
  \includegraphics[scale=1]{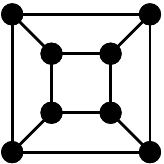}
}
\qquad
\subcaptionbox{}{
  \includegraphics[scale=1]{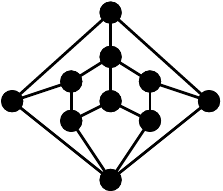}
}
\caption{Simple irreducible \acp{MUQ} of order 8 and 10}
\label{fig:simpleirreds}
\end{figure}

These two \acp{MUQ} are shown in \autoref{fig:simpleirreds}.
As there is also an irreducible with parallel edges depicted in \autoref{fig:parallelirred}, we have an obvious upper bound\cite{kapolnai2012periodica}:

\tetel{%
The order of the smallest irreducible \ac{MUQ} containing parallel edges is at most 14.
}

\section{Prerequisites}

In the proof of \autoref*{theo:nosmallparallel}, we use some concepts and properties of plain graphs.
Let us denote the number of vertices, edges and faces, by $n$, $e$ and $f$, respectively.
Euler's formula holds for any plane graph:
$n-e+f=2$.
Applied to a \ac{MUQ}, as every face has 4 boundary edges and every edge is counted twice, we have
\begin{equation}
\label{eq:Euler}
e%
=2n-4.
\end{equation}
Consequently, the sum of the degrees of a \ac{MUQ} of order $n$ equals to $4n-8$.

\autoref{eq:Euler} also implies
that the minimum degree of a
\ac{MUQ} is either 1, 2 or 3, because if it had only vertices of degree at least 4, then it would have at least $2n$ edges.

Let $G(V,E)$ denote a plane graph with vertex set $V$ %
 and edge set $E$, where $V$ contains points, $E$ contains arcs on the plane.
 The plane graph $G(V,E)$ is the \emph{embedded subgraph} of the plane graph $G'(V',E')$ if $V\subseteq V'$ and $E\subseteq E'$.

\section{Proof of \autoref*{theo:nosmallparallel}}

\begin{proof}
Let $G$ be an irreducible \ac{MUQ} with parallel edges such that it is of minimum order.
Suppose $G$ has $k$ parallel edges between vertices $v$ and $w$.
The $k$ parallel edges divide the surface into $k$ regions.
Let us select a region with the fewest vertices inside.
Now let us prepare the embedded subgraph $H$ of the \ac{MUQ} $G$ by selecting the subgraph spanned by the vertices from inside the selected region and the vertices $v$ and $w$, then removing $k-1$ parallel edges between $v$ and $w$.
For example, in \autoref{fig:parallelirred}, if the right hand side is $G$, one of the graphs on the left hand side could be isomorphic to $H$, and $k=2, x=v, y=w$.

\begin{figure*}
\centering
\includegraphics[]{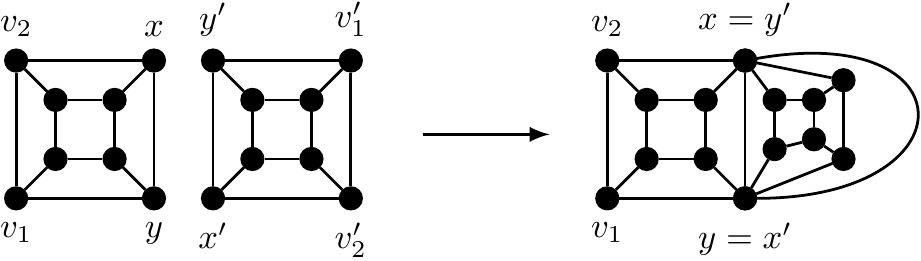} %
\caption[Irreducible quadrangulation with parallel edges]{Irreducible \ac{MUQ} with parallel edges created by repeating a simple one}
\label{fig:parallelirred}
\end{figure*}

It is easy to see that the plane graph $H$ has the following properties:
\begin{enumerate*}[label=(\roman{*})]
\item it is a \ac{MUQ},
\item $v$ and $w$ are still adjacent and there is exactly one edge between them,
\item \label{halfdeg} every other vertex than $v$ and $w$ has degree at least 3.
\end{enumerate*}
We call a plane graph satisfying these three conditions a \emph{half} with marked adjacent vertices $v$ and $w$ throughout this proof.
Let $d_K(z)$ denote the degree of vertex $z$ in graph $K$.
In addition to these three defining properties, it is easy to see that any \ac{MUQ} $Q$ also has property
\begin{enumerate*}[resume,label=(\roman{*})]
\item: \label{halfprop}
$d_Q(x)+d_Q(y)\geq 3$ for any two adjacent vertices $x$ and $y$.
\end{enumerate*}

An irreducible \ac{MUQ} with parallel edges can be built from a {half} $F$ with marked adjacent vertices $x$ and $y$, supposing there is only one edge between $x$ and $y$, as follows.
We clone $F$, rotate the clone by 180 degrees preserving its orientation, and stick together $F$ and the clone by unifying vertex $x$ and the clone image of $y$ denoted $y'$, unifying vertex $y$ and the clone image of $x$ denoted $x'$, and unifying the outer edge $xy$ and its clone image. %
So the unifications remove the clone vertices $x'$ and $y'$ and the clone of the edge $xy$.
Then we duplicate the original edge $xy$ outside in order to restore the quadrilateral property.
The process is illustrated in \autoref{fig:parallelirred}.
We need to show that the resulting graph denoted by $P$ is irreducible.
Clearly the order of $P$ is at least 4, so irreducibility is equivalent to having minimum degree of 3.
It is already a half, so every other vertices than $x$ or $y$ has degree at least 3.
After the unification, the vertex $x$ is now connected with the neighbours of $y'$ including $x'$ because of the additional parallel edge, so $d_P(x) = d_F(x)+d_F(y) \geq 3$ from property \ref{halfprop} of the halves.
Similarly, $d_P(y) = d_F(x)+d_F(y) \geq 3$.

Now we can prove that $H$ is simple.
By property \ref{halfdeg}, $H$ is not isomorphic to the \ac{MUQ} $P_2$.
Note that it is easy to see that if a half $H$ with marked adjacent vertices $v$ and $w$ has parallel edges, then it also contains a smaller half.
The contained half can be found using a similar method used to define the half $H$ inside the irreducible $G$, although instead of selecting the minimal region, any region can be selected which does not contain $v$ and $w$ (because $d_H(v)$ or $d_H(w)$ can be less than 3).
So suppose indirectly that $H$ has parallel edges.
Then there is a \emph{half} $F$ contained in $H$, obviously smaller than $H$.
So an irreducible \ac{MUQ} with parallel edges could be built from $F$ as described above, but smaller than $G$, contradicting to its minimality.

Now we are going to prove that $d_H(v)+d_H(w)\geq 5$.
The minimum degree of a simple quadrangulation is at least 2.
If $d_H(v)=d_H(w)=2$ holds in a simple quadrangulation for two adjacent vertices, than it can only be the circle $C_4$ of length 4.
However, by property \ref{halfdeg}, $H$ cannot be isomorphic to $C_4$.

If the half $H$ has $n$ vertices, than the sum of its degrees is $4n-8$, so we have
\begin{eqnarray*}
4n-8 = d_H(v)+d_H(w) + \sum_{z\neq v,w} d_H(z) \geq 5 + 3(n-2),
\end{eqnarray*}
implying $n\geq 7$.
So there are at least 5 vertices in the inside region of $G$, the same outside, plus $v$ and $w$, added up to 12.
\end{proof}

\section{Conclusions}

\autoref*{theo:nosmallparallel} means that when the software Plantri generates the simple irreducibles for $n<12$, then it also generates every irreducible.
Consequently, the cardinalities of the simple irreducibles of different orders published by Brinkmann et.\ al.\cite{Brinkmann2005} (see their Table 2 titled ``Simple quadrangulations with minimum degree 3'') does not exclude any parallel irreducible for $n<12$.

Up to the best knowledge of the authors, the identity of the smallest parallel irreducible is still unknown.
Although there is a primitive implementation to enumerate every \ac{MUQ}, it is practically unusable for $n\geq 12$ because of its very low efficiency.
A better way to find it could be extending Plantri to enumerate efficiently every \ac{MUQ} (not just the simple ones), using extended operations\cite{kapolnai2012periodica}.

\section*{Acknowledgements}

Being part of project Bioklíma, the research has been supported the Hungarian Government, managed by the National Development Agency, and financed by the Research and Technology Innovation Fund. 

The second author has been supported by OTKA grant 104601.

\newcommand{\ZsInitial}{Zs}
\bibliographystyle{abbrv}
\bibliography{suh,mypub}

\end{document}